\documentclass[11pt]{article}
\usepackage{amssymb, amsmath, amsthm}

\begin{document}
\parindent 0in
\parskip 1 em
\date{}
\title{\text{\bf{\normalsize{THE BOCHNER FORMULA VIA VOLUME VARIATIONS}}}}
\author{Christopher Lin}

\maketitle 

\begin{abstract}
 In this short paper, we re-derive the Bochner formula for the Laplacian by considering local variations of volume.  
 The derivation is rooted in the fact that the Laplacian of a function measures the volume 
variation along the flow of the gradient vector of the function.  Possible extensions of this approach/technique are also discussed.  While 
the value of this approach may be limited in terms of research, we think it definitely has pedagogical value.        
\end{abstract}

\section{Introduction}
\quad The classic Bochner formula 
\begin{equation}\label{bochner}
 \frac{1}{2}\Delta |\nabla f|^2 = \text{Rc}\,(\nabla f, \nabla f) + |\nabla^2 f|^2 + \langle\nabla f, \nabla\Delta f\rangle
\end{equation}
for any smooth function $f$ on a Riemannian manifold $M$ is one of the most important formulas in geometric analysis.  
In particular, it is a formula displaying  
the interplay between the geometry (the Ricci curvature) and the function theory (harmonic functions, eigenfunctions, etc) on a Riemannian manifold.  
The traditional derivation of (\ref{bochner}) is a direct computation of $\Delta |\nabla f|^2$, 
involving switching covariant derivatives and invoking curvature terms.  However straight-forward this traditional derivation may be, 
it seems to lack geometric insight.  We present here a derivation that is geometrically-inspired at the outset.  

\quad A motivation for this paper comes from the study of harmonic functions on a Riemannian manifold.  A harmonic function is a function $f$ satisfying 
Laplace's equation: $\Delta f = 0$.  The Laplacian $\Delta$ is defined as $\Delta = \text{div}\,\nabla$.  A geometric way of understanding
 Laplace's equation lies in the gradient $\nabla f$.  Let $\varphi_t$ be the local 
$1$-parameter group of diffeomorphisms generated by a vector field $X$ on $M$.  The divergence $\text{div}\,X$ of $X$ 
is a measure of how the volume of any small subset $\Omega$ on $M$ changes as it is carried through by $\varphi_t$, as expressed by the following elementary 
formula
\begin{equation}\label{var1}
 \frac{d}{dt}\int_\Omega \,dv = \int_\Omega \phi_t^* \big(\text{div}\,X\, dv\big)  .
\end{equation}
Let $f$ be a function on a Riemannian manifold $M$ with metric $g = \langle\,\,\,,\,\,\,\rangle$.  When $X = \nabla f$, we see that the Laplacian 
appears on the righ-hand side.  Therefore Laplace's equation actually says that small parcels on $M$ move without changing its volume along 
the gradient vector field $\nabla f$.  This immediately gives some geometric insight, without doing more analysis.  It is well-known that the 
curvature (in particular, Ricci curvature) of a Riemannian manifold is related to volume (in a local sense).  
Therefore, the curvature on a Riemannian manifold should somehow restrict the presence of harmonic functions with certain growth rates.  In fact, this is 
exactly the case, as vindicated by the many theorems in this direction.\footnote{For a survey of results and techniques, see for example, \cite{PLi2}.}

\quad  Equation (\ref{var1}) by itself does not involve any curvature at all.  What we need to do is to 
compute the second variation of volume.  In fact, our strategy is to rewrite the volume form as $dv = \sqrt{g}\,dx_1 \wedge\cdots\wedge dx_n$ in a 
specially chosen local coordinate system $(x_1,...,x_n)$ adapted to the gradient $\nabla f$, compute two variations of this form of $dv$ in the direction 
of the local diffeomorphism $\varphi_t$, and then compare (equate) it to the same computations done in the flavor of (\ref{var1}).  The 
Bochner formula (\ref{bochner}) then pops right away.  Aligning with our 
point of view, (\ref{bochner}) is a reflection that the Ricci curvature controls the acceleration\footnote{as opposed to the velocity, or 
simply the rate of change 
of volume, which does not depend on curvature} 
of the volume of a small parcel on $M$ as it is carried along 
the flow of $\nabla f$.      

\quad The method outlined above may seem  
like a tautological redundency to simply re-derive the variations in terms of the divergence by using local coordinates.  However, to compute the variations 
of volume as in (\ref{var1}) is more of an exercise in Lie derivatives, and can be done quite easily.  On the other hand, 
computing the variations in terms of the specially chosen local coordinate system will be more difficult and is really quite different.  In fact, the 
all-important appearance of the Ricci curvature in (\ref{bochner}) comes from the computation in local coordinates.  This part of the computation 
is reminiscent of the standard computation of the first and second variations of area along geodesics.\footnote{for example, see \cite{PLi}}  
However, here our choice of local coordinate system is specific to the function $f$ and although we encounter the same terms obtained in the standard 
computation of the variations of area, we also have to further expand out some of these terms in more detail.      

\quad In what follows we will demonstrate the computations in the derivation outlined above.  We want to emphasize that besides the 
scrupulous choice of a local coordinate system, the rest of our computations does not involve any tricks - only some patience is required.  The most 
important message here, of course, is the fact that the Bochner formula (\ref{bochner}) is actually just a statement about how volume changes 
(accelerates) along
 the flow 
of $\nabla f$.  One can view this as a justification, on an intuitive level, of the many results in geometric analysis.  This method can also be 
applied to deriving similar formulas for other operators related to the Laplacian.  In fact, at the end of the paper we will discuss  
a similar Bochner formula for the so-called drift Laplacian in light of the method presented here. 

\text{\bf{Acknowlegement}}\, The author would like to acknowlege the invaluable lecture notes of Professor Peter Li, and whose style in 
geometric analysis has continued to inspire over the years.   

\section{Notations and Basic Facts}
\quad In this section we collect the notations and basic facts that will be used in our computations.  For a more comprehensive review 
of Riemannian geometry, please see for example, 
the self-contained lecture notes \cite{PLi} of P. Li or the book \cite{Jost} by Jost.  
  
\quad Let $M$ be a Riemannian manifold of dimension $n$ with metric $g = \langle\,\,\,,\,\,\,\rangle$.  Unless mentioned otherwise, 
for the computations in this 
paper we will refrain from the usual Einstein summation convention, 
since in our later computations we wish to distinguish between coordinates on a hypersurface and coordinates 
in an open set of the manifold $M$.    
The Levi-Civita connection of $g$ is denoted by $\nabla$, which acts on 
vector fields on $M$.  The Curvature tensor of the Levi-Civita conection is defined by 
\begin{equation}\label{curvature}
 R(X,Y)Z = \nabla_X\nabla_Y Z - \nabla_Y\nabla_X Z - \nabla_{[X,Y]}Z.
\end{equation}
The Ricci curvature is a symmetric $2$-tensor defined by 
\begin{equation}
 \text{Rc}\,(X,Y) = \,\text{tr}\langle R(\,\,,X)Y,\,\,\rangle = \sum_{i=1}^n\langle R(e_i, X)Y,e_i\rangle , 
\end{equation}
where $\{e_i\}$ is any point-wise orthonomal frame.
The divergence of a vector field $X$ is the function defined by 
\begin{equation}
 \text{div}\,X = \text{tr}(Y \mapsto \nabla_YX) = \sum_{i=1}^n \langle\nabla_{e_i} X, e_i\rangle.
\end{equation}

\quad The gradient of a function $f$ on $M$ is also denoted by $\nabla$.  The Hessian of $f$ is the symmetric $2$-tensor defined by 
\begin{equation}
 \nabla^2 f \,(X,Y) = \langle \nabla_X \nabla f\,,\, Y\rangle.
\end{equation}

\quad The Lie derivative depends only on the differential structure of the manifold $M$, and is defined on a differential form $\omega$ by 
\[
 L_X\omega = \frac{d}{dt}\varphi_t^*(\omega)\Big|_{t=0}
\]
at any point on $M$, where $X$ is a vector field on $M$ with $\varphi_t$ its associated local $1$-parameter group of diffeomorphisms.  The Lie derivative 
on differential forms enjoys Cartan's magic formula
\begin{equation}\label{Cartan}
 L_X\omega = d(X\lrcorner\omega) + X\lrcorner d\omega.
\end{equation}

\section{The Derivation} Although it is a standard computation, we will first recall the derivation of (\ref{var1}) for the sake of completeness.  
If we don't evaluate at ${t=0}$, we have
\begin{equation}\label{minor}
 \frac{\partial}{\partial t}\varphi_t^*(dv) = \varphi_t^*\big(L_X dv\big). 
\end{equation}        
By (\ref{Cartan}) and the fact that $dv$ is an $n$-form, we have  
$L_X dv = d(X\lrcorner dv) + X\lrcorner d(dv) = d(X\lrcorner dv)$.  We can then write $L_X dv = F\,dv$ for some function $F$ on $M$.  
At any point $p\in M$ with respect to any orthonormal frame, we also have
 $\{e_1,...,e_n\}$, we have $F = d(X\lrcorner dv)(e_1,...,e_n)$.  Instead of doing the calculation by brute force, 
we may choose a geodesic frame at $p$, i.e. $\langle e_i, e_j\rangle = \delta_{ij}$ and $\nabla_{e_i}e_j = 0$ at $p$ for all $i$ and $j$.  Then at $p$ we 
have 
\begin{align}
 F &= d(X\lrcorner dv)(e_1,...,e_n) \notag\\
  &= e_i\big(X\lrcorner dv(e_1,...,X,...,e_n)\big)  \hskip 1cm (\text{$X$ placed at the $i$th slot}) \notag\\
  &= e_i\big(a_j\,dv(e_1,...,e_n)\big)               \hskip 1.2cm (\text{writing $X=a_j e_j$}) \notag\\
  &= e_i(a_i) = \text{div}\, X. \notag
\end{align}
Substituting this into (\ref{minor}) and integrate, we have (\ref{var1}).  By differentiating again in the same way, we have 
\begin{align}\label{leftside}
 \frac{\partial^2}{\partial t^2}\varphi_t^*(dv)\Big|_{t=0} &= L_X(\text{div}\,X\,dv) \notag\\
  &= X\text{div}\,X \,dv + \text{div}X\,L_X(dv) \notag\\
  &= \big(X\text{div}\,X + (\text{div}\,X)^2\big)\,dv. 
\end{align}

\quad Next we will proceed with the computation of $\frac{\partial^2}{\partial t^2}\varphi_t^*(dv)\Big|_{t=0}$ in local coordinates.  First we assume 
that $|\nabla f|\ne 0$ at a point $p\in M$.  Since $f$ is smooth, there exists an open neighborhood $U$ containing $p$ over which $|\nabla f|\ne 0$.  
Moreover, we may assume that $U$ is foliated by level sets of $f$.  We want a coordinate system on $U$ adapted to this foliation, and whose 
ensuing metric tensor is easy to handle.  To this end, we consider 
the vector field $X = \frac{\nabla f}{|\nabla f|^2}$.

\newtheorem{lemma}{Lemma}
\begin{lemma}\label{sametime}
 Consider any two level sets $\{f=C_1\}$ and $\{f=C_2\}$ in $U$.  The time it takes to go 
 from $\{f=C_1\}$ to $\{f=C_2\}$ along integral curves of $X$ is the 
 same for every point on $\{f=C_1\}$. 
\end{lemma}

\text{\bf{Proof}.} \, Let $\gamma$ be an integral curve of $X$, i.e. $\dot{\gamma} = X$, starting at a point $p\in\{f=C_1\}$ and 
ending at a point $q\in\{f=C_2\}$.  Let $T$ be the time 
it takes to reach $\{f=C_2\}$.  Then we see that 
\begin{align}
 C_2 - C_1 = f(q) - f(p) = \int_\gamma df &= \int_0^T df(\dot{\gamma})\,dt \notag\\
  &= \int_0^T df(\nabla f / |\nabla f|^2) \,dt \notag\\
  &= \int_0^T \langle\nabla f, \frac{\nabla f}{|\nabla f|^2} \rangle\, dt \notag\\
  &= \int_0^T dt = T.\notag
\end{align}  
Since $p\in\{f=C_1\}$ is arbitrary, we see that the time is always the difference $C_2 - C_1$. \qed

With the Lemma above, we can now define our local coordinate system as follows.  Consider the original point $p$, and let $\{f= C_0\}$ be the level set 
containing $p$.  Since $\{f= C_0\} \cap U$ is a hypersurface in $M$, we can consider a local coordinate system $(x_2,...,x_n)$ on $\{f= C_0\}$.  We can 
assume assume that $U$ is small enough so that $\{f= C_0\} \cap U$ is covered entirely by such a local coordinate system.  Then for every $q\in U$ there is a 
unique point $\tilde{q}\in \{f= C_0\} \cap U$ from which the integral curve of $X$ reaches $q$.  We now define the coordinate system $(x_1,x_2,...,x_n)$ on 
$U$, where at any $q\in U$ the coordinate $x_1$ is the value of $f$ and $x_2,...,x_n$ are the inherited coordinates from that of $\tilde{q}$.  
By Lemma \ref{sametime}, we thus have $\partial_1 = \frac{\partial}{\partial x_1} = X$.  Since 
gradient vectors are orthogonal to level sets and $|X|^2 = |\nabla f|^{-2}$, the metric is given by 
\begin{equation}\label{Gij}
 G_{ij} = 
 \begin{cases}
  |X|^2 \hskip 1cm i=j=1 \\
  0 \hskip 1.6cm i\ne 1 \,\text{and}\, j=1, \,\text{or}\, j\ne 1 \,\text{and}\, i=1 \\
  g_{ij} \hskip 1.5cm i,j \ne 1
 \end{cases}
\end{equation}
over $U$, and $g_{ij}$ denotes the metric on the level sets of $f$ in $U$.  Therefore the inverse is given by 
\begin{equation}\label{Gupperij}
 G^{\,ij} = 
 \begin{cases}
  |X|^{-2} \hskip 1cm i=j=1 \\
  0 \hskip 1.6cm i\ne 1 \,\text{and}\, j=1, \,\text{or}\, j\ne 1 \,\text{and}\, i=1 \\
  g^{ij} \hskip 1.5cm i,j \ne 1
 \end{cases}
\end{equation} 
over $U$.  We can further assume that the coordinate vectors $\{\partial_2,...,\partial_n\}$ are orthonormal at $p$.\footnote{For example we can choose 
$(x_2,...,x_n)$ to be normal coordinates on $\{f= C_0\} \cap U$ centered at $p$.}  This will be useful in the computations below.  

\quad From the choice of the local coordinate system above, we see that $\varphi_t^* (dx_i) = d(x_i\circ\varphi_t) = dx_i$ for all $i=1,...,n$. 
We also write $\varphi_t^*(G)|_p = G(p,t)$, where $G(p,0) = G(p)$ and $G$ here denotes the determinant of the metric $(G_{ij})$.  Similar notations will be used for the metric $g_{ij}$.  Then by also 
writing $J(p,t) = \sqrt{G(p,t)}/\sqrt{G(p)}$, at $p$ we can write 
\begin{align}
 \varphi_t^*(dv) &= \sqrt{G(p,t)}\,dx_1\wedge\cdots dx_n \notag\\
  &= J(p,t)\,dv. \notag
\end{align}
Recalling Jacobi's formula $\frac{d}{dt}(\det A) = \det A \,\text{tr}(A^{-1}\frac{dA}{dt})$ for an invertible matrix $A$, we have
\begin{equation}\label{1stpar}
 \frac{\partial}{\partial t}J(p,t) = \frac{1}{2}J(p,t)G^{ij}\frac{\partial G_{ij}}{\partial t}, 
\end{equation}    
 which implies that 
\begin{equation*}
 \frac{\partial}{\partial t}J(p,t)\Big|_{t=0} = \frac{1}{2}G^{ij}(p)\frac{\partial G_{ij}}{\partial t}\Big|_{t=0}.
\end{equation*}
From (\ref{Gij}), (\ref{Gupperij}), and (\ref{1stpar}), we then have\footnote{Note that it is unambiguous to use the Einstein 
summation convention for terms involving $G_{ij}$ and $G^{ij}$ below, and we also use the well-known formula 
$\frac{\partial G^{ij}}{\partial t} = -G^{ik}\frac{\partial G_{kl}}{\partial t} G^{lj}$.}
 (we will omit the evaluation at $t=0$ on the righ-hand side below)
\begin{align}\label{2ndpar1}
 \frac{\partial^2 J}{\partial t^2}\Big|_{t=0} &= \frac{1}{2}\Big(\frac{\partial J}{\partial t}G^{ij}\frac{\partial G_{ij}}{\partial t} + 
 \frac{\partial G^{ij}}{\partial t}\frac{\partial G_{ij}}{\partial t} + G^{ij}\frac{\partial^2G_{ij}}{\partial t^2}\Big) \notag\\
 &= \frac{1}{4}\Big(|X|^{-2}\frac{\partial |X|^2}{\partial t}+\sum_{i=2}^{n}\frac{\partial g_{ii}}{\partial t}\Big)^2 - 
 \frac{1}{2}|X|^{-4}\Big(\frac{\partial |X|^2}{\partial t}\Big)^2 - \frac{1}{2}\sum_{i,j = 2}^n \Big(\frac{\partial g_{ij}}{\partial t}\Big)^2 \notag\\
 &\,\,\,+ \frac{1}{2}|X|^{-2}\frac{\partial^2|X|^2}{\partial t^2} + \frac{1}{2}\sum_{i=2}^n\frac{\partial^2 g_{ii}}{\partial t^2}.
\end{align}

An elementary identity that will be used repeatedly is the following.
\begin{lemma}\label{useful}
 $Y(|X|^2) = -2|X|^2 \,\nabla^2 f(Y,X)$,  for any vector $Y$,
\end{lemma}

\text{\bf{Proof}.} \, We compute and see that 
\begin{align}
 Y(|X|^2) &= 2\langle\nabla_Y X, X\rangle \notag\\
  &= 2\langle\nabla_Y (|X|^2\nabla f), X \rangle \notag\\
  &= 2\langle Y(|X|^2)\nabla f + |X|^2\nabla_Y\nabla f, X\rangle\notag\\
  &= 2Y(|X|^2) + 2|X|^2\langle\nabla_Y\nabla f,X\rangle. \notag
\end{align}
Solving for $Y(|X|^2)$, the lemma follows. \qed

\quad Continuing, we note that 
\begin{align}
 \sum_{i=2}^n\frac{\partial g_{ii}}{\partial t}\Big|_{t=0} &= \sum_{i=2}^n X\langle\partial_i,\partial_i\rangle \notag\\ 
  &= 2\langle\nabla_X\partial_i,\partial_i\rangle \notag\\
  &= 2\langle\nabla_{\partial_i}X,\partial_i\rangle = 2\big(\text{div}\,X - |X|^{-2}\langle\nabla_X X,X\rangle\big),\notag
\end{align}
which implies 
\begin{align}\label{square}
 \frac{1}{4}\Big(|X|^{-2}\frac{\partial |X|^2}{\partial t}+\sum_{i=2}^{n}\frac{\partial g_{ii}}{\partial t}\Big)^2 &= 
 \frac{1}{4}\Big(|X|^{-2}X\langle X,X\rangle + 2\text{div}\,X - 2|X|^{-2}\langle \nabla_X X, X\rangle\Big)^2\notag\\
 &= \big(\text{div}\,X\big)^2.
\end{align}
We also have
\begin{align}\label{d2gii}
 \sum_{i=2}^n\frac{\partial^2 g_{ii}}{\partial t^2}\Big|_{t=0} &= 2 \sum_{i=2}^n X\langle\nabla_{\partial_i}X,\partial_i\rangle\notag\\
  &= 2\sum_{i=2}^n\langle\nabla_X\nabla_{\partial_i}X,\partial_i\rangle + 2\sum_{i=2}^n\langle\nabla_{\partial_i}X, \nabla_X\partial_i\rangle \notag\\
  &= 2\sum_{i=2}^n\langle R(X,\partial_i)X + \nabla_{\partial_i}\nabla_X X, \partial_i\rangle + 2\sum_{i=2}^n|\nabla_{\partial_i}X|^2 \notag\\
  &= -2\,\text{Rc}(X,X) + 2\,\text{div}\,\big(\nabla_X X\big) - 2|X|^{-2}\langle\nabla_X\nabla_X X, X\rangle + 2\sum_{i=2}^n|\nabla_{\partial_i}X|^2.
\end{align} 
Moreover,
\begin{align}\label{dgij}
 \frac{\partial g_{ij}}{\partial t}\Big|_{t=0} = X\langle\partial_i,\partial_j\rangle 
 &= \langle\nabla_X \partial_i,\partial_j\rangle + \langle\nabla_X\partial_j,\partial_i\rangle \notag\\
 &= \langle\nabla_{\partial_i} X, \partial_j\rangle + \langle\nabla_{\partial_j}X,\partial_i\rangle \notag\\
 &= \langle\nabla_{\partial_i}(|X|^2\nabla f),\partial_j\rangle + \langle\nabla_{\partial_j}(|X|^2\nabla f),\partial_i\rangle \notag\\
 &= |X|^2\langle\nabla_{\partial_i}\nabla f,\partial_j\rangle + |X|^2\langle\nabla_{\partial_j}\nabla f,\partial_i\rangle \notag\\
 &= 2|X|^2\,\nabla^2 f(\partial_i,\partial_j)
\end{align}
and
\begin{equation}\label{neweqn}
 \frac{\partial^2|X|^2}{\partial t^2} = XX|X|^2 = 2\langle\nabla_X\nabla_X X,X\rangle + 2|\nabla_X X|^2. 
\end{equation}

\quad On the other hand, via Lemma \ref{useful} we also have 
\begin{align}\label{leftpart}
 X\text{div}\,X &= X\big(\text{div}\,(|X|^2\nabla f)\big) \notag\\
  &= X\big(\langle \nabla|X|^2,\nabla f\rangle + |X|^2\Delta f\big) \notag\\
  &= X\big(\nabla f\langle X,X\rangle + |X|^2\Delta f\big) \notag\\
  &= X\big(|X|^{-2} X(|X|^2)\big) + X\big(|X|^2\Delta f\big) \notag\\
  &= X\big(|X|^{-2})X(|X|^2) + |X|^{-2} XX(|X|^2) + X(|X|^2)\Delta f + |X|^2 X\Delta f \notag\\
  &= -|X|^{-4}X(|X|^2)^2 + |X|^{-2}XX(|X|^2) + X(|X|^2)\Delta f + |X|^4 \langle\nabla f, \nabla\Delta f\rangle. \notag\\
\end{align}

In view of (\ref{2ndpar1}), (\ref{square}), (\ref{d2gii}),(\ref{dgij}), (\ref{neweqn}), and (\ref{leftpart}), 
by equating $\frac{\partial^2 J}{\partial t^2}\big|_{t=0}\,dv$ with (\ref{leftside}) we obtain 
\begin{align}\label{step1}
 -\frac{1}{2}|X|^{-4}&X(|X|^2)^2 - 2|X|^4\sum_{i,j=2}^n\nabla^2 f(\partial_i,\partial_j)^2 - \text{Rc}(X,X) + \text{div}\big(\nabla_X X\big) 
  + \sum_{i=1}^n|\nabla_{e_i}X|^2 \notag\\
  &= -|X|^{-4}X(|X|^2)^2 + |X|^{-2}XX(|X|^2) + X(|X|^2)\Delta f + |X|^4\langle\nabla f,\nabla\Delta f\rangle, 
\end{align} 
where for each $i=1,...,n$, $e_i$ is the unit vector in the direction of the coordinate vectors.  By Lemma \ref{useful} we have $X(|X|^2) = -2|X|^2\,\nabla^2 f(X,X)$, which again via Lemma \ref{useful} implies 
\[
XX(|X|^2) = 4|X|^2\,\nabla^2 f(X,X)^2 - 2|X|^2X\nabla^2 f(X,X).
\]
Substituting these two equations into (\ref{step1}) and after rearranging, we have 
\begin{align}\label{step2}
 \text{Rc}(X,X) + |X|^4\langle\nabla f,\nabla\Delta f\rangle 
  &= -2\nabla^2 f(X,X)^2 + 2X\nabla^2 f(X,X) + 2|X|^2\,\nabla^2 f(X,X)\Delta f \notag\\
   &\hskip 0.5cm - 2|X|^4\sum_{i,j=2}^n\nabla^2 f(\partial_i,\partial_j)^2 + \text{div}\big(\nabla_X X\big) + \sum_{i=1}^n|\nabla_{e_i}X|^2. 
\end{align} 

The proof of the Bochner formula will follow readily after we establish the following results.

\begin{lemma}\label{thru1}
 $\sum_{i=1}^n|\nabla_{e_i}X|^2 = |X|^4|\nabla^2 f|^2$
\end{lemma}

\text{\bf{Proof.}}\,
 \begin{align}
 \sum_{i=1}^n|\nabla_{e_i}X|^2 &= \sum_{i=1}^n\langle\,\sum_{j=1}^n\langle\nabla_{e_i}X, e_j\rangle e_j,\nabla_{e_i} X\rangle \notag\\
  &= \sum_{i,j=1}^n \langle\nabla_{e_i}X, e_j\rangle^2\notag\\
  &= \sum_{i,j=1}^n \Big(e_i(|X|^2)\langle\nabla f,e_j\rangle + |X|^2\langle\nabla_{e_i}\nabla f,e_j\rangle\Big)^2 \notag\\
  &= |X|^{-2}\sum_{i=1}^n e_i(|X|^2)^2 + 2\sum_{i=1}^n e_i(|X|^2)\nabla^2 f(e_i,X) + |X|^4|\nabla^2 f|^2 \notag\\
  &= |X|^4|\nabla^2 f|^2,\notag
\end{align}
where we have applied Lemma \ref{useful} with $Y=e_i$.  \qed

\begin{lemma}\label{thru2}
 \begin{align}
 \text{div}\big(\nabla_X X\big) = &-2X\nabla^2 f(X,X) + 4\nabla^2 f(X,X)^2 
 -2|X|^2\nabla^2 f(X,X)\Delta f \notag\\
  &-4|X|^2\sum_{i=1}^n\nabla^2 f(e_i,X)^2  +\frac{|X|^4}{2}\Delta|\nabla f|^2 \notag
 \end{align}
\end{lemma}

\text{\bf{Proof.}}\, First we note that 
\begin{align}
 \nabla_X X &= X(|X|^2)\nabla f + |X|^2\nabla_X\nabla f \notag\\
  &= -2|X|^2\,\nabla^2 f(X,X)\,\nabla f + |X|^2\nabla_X\nabla f, \notag
\end{align}
which implies 
\begin{align}\label{big}
 \text{div}\big(\nabla_X X\big) &= -2\nabla f \big(|X|^2\,\nabla^2 f(X,X)\big) -2|X|^2\nabla^2 f(X,X)\Delta f \notag\\ 
  &\hskip 0.5cm + \langle\nabla|X|^2,\nabla_X \nabla f\rangle + |X|^2\text{div}\big(\nabla_X\nabla f\big) \notag\\
  &= -2|X|^2\nabla f\big(\nabla^2 f(X,X)\big) -2\nabla f(|X|^2)\nabla^2 f(X,X) - 2|X|^2\nabla^2 f(X,X) \Delta f \notag\\
  &\hskip 0.5cm + \sum_{i=1}^n e_i(|X|^2)\langle e_i, \nabla_X\nabla f\rangle + |X|^2\text{div}\big(\nabla_X\nabla f\big) \notag\\
  &= -2X\nabla^2 f(X,X) + 4\nabla^2 f (X,X)^2 - 2|X|^2\nabla^2 f(X,X) \Delta f \notag\\
  &\hskip 0.5cm  -2|X|^2\sum_{i=1}^n\nabla^2 f(e_i,X)^2 + |X|^2\text{div}\big(\nabla_X\nabla f\big), 
\end{align} 
where we have applied Lemma \ref{useful} two times.  We further note that 
\begin{align}
 \nabla_X\nabla f = \sum_{i=1}^n\langle\nabla_X\nabla f, e_i\rangle \, e_i 
  &= |X|^2\sum_{i=1}^n\nabla^2 f(\nabla f, e_i)\,e_i \notag\\
  &= |X|^2\sum_{i=1}^n\langle\nabla_{e_i}\nabla f, \nabla f\rangle\,e_i \notag\\
  &= \frac{1}{2}|X|^2 e_i(|\nabla f|^2)\,e_i 
  = \frac{1}{2}|X|^2 \nabla |\nabla f|^2 ,\notag
\end{align}
which implies 
\begin{align} \label{support}
 \text{div}\big(\nabla_X\nabla f\big) &= \frac{1}{2}\langle\nabla|X|^2,\nabla|\nabla f|^2\rangle + \frac{|X|^2}{2}\Delta|\nabla f|^2 \notag\\
  &= \frac{1}{2}\langle e_i(|X|^2)e_i\,,\,2\nabla^2 f(\nabla f, e_j)e_j\rangle+ \frac{|X|^2}{2}\Delta|\nabla f|^2\notag\\
  &= \sum_{i=1}^n e_i(|X|^2)\nabla^2 f(\nabla f,e_i) + \frac{|X|^2}{2}\Delta|\nabla f|^2\notag\\
  &= -2\sum_{i=1}^n \nabla^2 f(e_i,X)^2+ \frac{|X|^2}{2}\Delta|\nabla f|^2.
\end{align}
Substituting (\ref{support}) into (\ref{big}), the Lemma follows. \qed 

\quad Finally, applying Lemma \ref{thru1} and Lemma \ref{thru2} in (\ref{step2}) we see that
\begin{align}\label{finally}
 \text{Rc}(X,X) + |X|^4\langle\nabla f,\nabla\Delta f\rangle 
  &= -2\nabla^2 f(X,X)^2 + 2X\nabla^2 f(X,X) + 2|X|^2\,\nabla^2 f(X,X)\Delta f \notag\\
   &\hskip 0.5cm - 2|X|^4\sum_{i,j=2}^n\nabla^2 f(\partial_i,\partial_j)^2 -2X\nabla^2 f(X,X) + 4\nabla^2 f(X,X)^2 \notag\\ 
   &\hskip 0.5cm -2|X|^2\nabla^2 f(X,X)\Delta f -4|X|^2\sum_{i=1}^n\nabla^2 f(e_i,X)^2  +\frac{|X|^4}{2}\Delta|\nabla f|^2 \notag\\
   &\hskip 0.5cm + |X|^4|\nabla^2 f|^2 \notag\\
   &= 2\nabla^2 f(X,X)^2 - 2|X|^4\sum_{i,j=2}^n\nabla^2 f(\partial_i,\partial_j)^2 \notag\\
    &\hskip 0.5cm -4|X|^2\Big(\nabla^2 f(e_1,X)^2 + \sum_{i=2}^n\nabla^2 f(\partial_i,X)^2\Big)\notag\\
     &\hskip 0.5cm + \frac{|X|^4}{2}\Delta|\nabla f|^2+ |X|^4|\nabla^2 f|^2 \notag\\
     &= -2|X|^4\sum_{i,j=2}^n\nabla^2 f(\partial_i,\partial_j)^2 -2|X|^4\nabla^2 f(e_1,e_1)^2 \notag\\
      &\hskip 0.5cm - 4|X|^4\sum_{i=2}^n\nabla^2 f(\partial_i, e_1)^2 + \frac{|X|^4}{2}\Delta|\nabla f|^2+ |X|^4|\nabla^2 f|^2 \notag\\
      &= \frac{|X|^4}{2}\Delta|\nabla f|^2 -  |X|^4|\nabla^2 f|^2.
\end{align}
Under our original assumption that $|\nabla f|\ne 0$ at $p$, dividing both sides of (\ref{finally}) by $|X|^4$ 
immediately yields the Bochner formula (\ref{bochner}).  If $\nabla f =0$ in an open neighborhood 
containing $p$, then (\ref{bochner}) holds trivially at $p$.  If $p$ is a critical point of $f$ that is 
not contained in an open neighborhood of critical points, formula (\ref{bochner}) still holds at 
$p$ by passing through the limit approaching $p$, since all objects are smooth.  Our derivation of the Bochner formula is now complete.    

\section{Further Comments}
\quad The technique demonstrated above is readily applicable to the derivation of other Bochner formulas.  The drift Laplacian on a 
Riemannian manifold $M$ is defined by 
\[
 \Delta_\phi = \Delta - \langle\nabla\phi, \nabla\,\,\,\rangle 
\]
for any function $\phi$ on $M$.  Clearly, when $\phi \equiv 0$ it reduces to the usual Laplacian.  On a compact manifold, 
$-\Delta_{\phi}$ is a non-negative definite elliptic 
operator that is (formally) self-adjoint with respect to the measure $e^{-\phi}\,dv_g$, i.e.
\[
   \int_M f(\Delta_\phi g)\,e^{-\phi}\,dv = \int_M g(\Delta_\phi f)\, e^{-\phi}\,dv 
\]  
for all $f,g\in C^{\infty}(M)$.  Generalizing the Laplacian case, we observe that for any vector field $X$ on $M$, we have 
\begin{equation}\label{vardrift}
 \frac{\partial}{\partial t}\big(e^{-\phi}\,dv\big)\Big|_{t=0} = \Big(-\langle X,\nabla\phi\rangle + \text{div}\,X\Big)e^{-\phi}\,dv,
\end{equation}
along the local $1$-parameter group of diffeomorphisms $\phi_t$ generated by $X$.  When $X = \nabla f$, the term in the paranthesis  of (\ref{vardrift}) 
becomes $\Delta_\phi f$.  Therefore, the drift Laplacian of a function $f$ also 
measures the rate of change of the $\phi$-volume $\int_\omega\, e^{-\phi}\,dv$ of any small parcel on $M$, along the flow of the vector field $\nabla f$.  
In particular, the drifted Laplace's equation $\Delta_\phi f = 0$ says that the $\phi$-volume change is zero along the flow of $\nabla f$.  By 
differentiating again with respect to $t$ and using the  
formulas in the last section (with the same coordinate system adapted to the level sets of $f$, where $X = \nabla f/|\nabla f|^2 = \partial_1$), one can readily derive the more general Bochner formula for the drift Laplacian:
\begin{equation}\label{driftbochner}
 \frac{1}{2}\Delta_\phi |\nabla f|^2 = |\nabla^2 f|^2 + \langle\nabla f, \nabla\Delta_\phi f\rangle + \text{Rc}_\phi\,(\nabla f,\nabla f) , 
\end{equation}    
where $\text{Rc}_\phi = \text{Rc} + \nabla^2\phi$ is the so-called Bakry-\'{E}mery tensor.\footnote{or more precisely, the $\infty$-Bakry-\'{E}mery tensor, 
first introduced in \cite{BE}}  
It is now well-known\footnote{For example, see \cite{WW}.} 
that the Bakry-\'{E}mery tensor is related to the $\phi$-volume in a more subtle, but analogous way as the Ricci tensor is related to the usual 
Riemannian volume.  Therefore, 
(\ref{driftbochner}) simply reflects the drift Laplacian's characterization of the rate of change of $\phi$-volume.  In particular, the 
Bakry-\'{E}mery tensor also controls the acceleration of $\phi$-volume along the flow of $\nabla f$.

\quad It would be interesting to find other examples of operators to which we can apply the same technique above and derive associated Bochner 
formulas.  More precisely, one would need to find other notions of volume on manifolds which can be related to operators.  Although this technique may not 
necessarily give us new formulas, it would definitely provide a more geometric interpretation to existing ones.

\vskip 0.8cm

    \text{\scshape{Department of Mathematics,
             Case Western Reserve University}},\\
            \text{\scshape{10900 Euclid Avenue-Yost Hall Room 220,
             Cleveland, OH}} 44106-7058\\
\text{\it{E-mail Address}}: ccl37@case.edu

\end{document}